\documentclass{gtart_a}
\pdfoutput=1

\usepackage[all]{xy}

\title{Rational maps and string topology}

\author{Sadok Kallel}
\givenname{Sadok}
\surname{Kallel}
\address{Laboratoire Painleve\\Universite de Lille I\\\newline
Villeneuve d'Ascq\\France}
\email{sadok.kallel@math.univ-lille1.fr}

\author{Paolo Salvatore}
\givenname{Paolo}
\surname{Salvatore}
\address{Dipartimento di matematica\\
Universita di Roma ``Tor Vergata''\\\newline
Roma\\Italy}
\email{salvator@mat.uniroma2.it}

\volumenumber{10}
\issuenumber{}
\publicationyear{2006}
\papernumber{39}
\lognumber{0363}
\startpage{1579}
\endpage{1606}

\doi{}
\MR{}
\Zbl{}

\keyword{mapping space}
\keyword{rational map}
\keyword{string product}
\subject{primary}{msc2000}{58D15}
\subject{secondary}{msc2000}{55R20}
\subject{secondary}{msc2000}{26C15}

\received{23 September 2003}
\revised{28 August 2006}
\accepted{11 September 2006}
\proposed{Ralph Cohen}
\seconded{Haynes Miller, Bill Dwyer}
\published{28 October 2006}
\publishedonline{28 October 2006}
\corresponding{}
\editor{CPR}
\version{}

\arxivreference{math.AT/0309137}

\let\xysavmatrix\xymatrix
\def\xymatrix{\disablesubscriptcorrection\xysavmatrix}
\AtBeginDocument{\let\bar\wbar\let\tilde\wtilde}
\numberwithin{equation}{section}
\makeop{EXP}
\makeop{exp}
\makeop{Im}
\makeop{Homeo}

\makeatletter
\def\cnewtheorem#1[#2]#3{\newtheorem{#1}{#3}[section]
\expandafter\let\csname c@#1\endcsname\c@theorem}

\newtheorem{theorem}{Theorem}[section]
\newtheorem{proposition}[theorem]{Proposition}
\newtheorem{lemma}[theorem]{Lemma}
\newtheorem{corollary}[theorem]{Corollary}
\theoremstyle{definition}
\newtheorem{remark}[theorem]{Remark}
\newtheorem{example}[theorem]{Example}

\newtheorem{definition}[theorem]{Definition}
\makeatother  

\newcommand{\bth}{\begin{theorem}}
\newcommand{\bpr}{\begin{proposition}}
\newcommand{\epr}{\end{proposition}}
\newcommand{\bco}{\begin{corollary}}
\newcommand{\eco}{\end{corollary}}
\newcommand{\bde}{\begin{definition}}
\newcommand{\ede}{\end{definition}}
\newcommand{\ble}{\begin{lemma}}
\newcommand{\ele}{\end{lemma}}
\newcommand{\bre}{\begin{remark}\rm}
\newcommand{\ere}{\end{remark}}
\newcommand{\bex}{\begin{example}\rm}
\newcommand{\eex}{\end{example}}

\def\la#1{\hbox to #1pc{\leftarrowfill}}
\def\ra#1{\hbox to #1pc{\rightarrowfill}}
\def\fract#1#2{\raise4pt\hbox{$ #1 \atop #2 $}}

\def\za{\qed}
\def\lrar{{\ra 2}}

\def\map#1{\mathrm{Map}_{#1}}
\def\bmap#1{\mathrm{Map}^*_{#1}}
\def\rat#1{\mathrm{Rat}_{#1}}
\def\hol#1{\mathrm{Hol}_{#1}}

\def\map#1{\mathrm{Map}_{#1}}

\def\tensor{\otimes}

\def\bbz{{\mathbb Z}}
\def\bbf{{\mathbb F}}
\def\bbp{{\mathbb P}}
\def\bbr{{\mathbb R}}

\def\bbc{{\mathbb C}}
\def\bbr{{\mathbb R}}
\def\bbq{{\mathbb Q}}

\def\bbh{{\mathbb H}}
\def\bbe{{\mathbb E}}

\def \mp{{\mathfrak p}}
\def \mg{{\mathfrak g}}
\def \mr{{\mathfrak r}}
\def \mn{{\mathfrak n}}
\def \mP{{\mathfrak P}^-}

\def\ot{\otimes}

\begin{document}

\begin{abstract}
We apply a version of the Chas--Sullivan--Cohen--Jones product on the
higher loop homology of a manifold in order to compute the homology of
the spaces of continuous and holomorphic maps of the Riemann sphere
into a complex projective space.  This product makes sense on the
homology of maps from a co--$H$ space to a manifold, and comes from a
ring spectrum.  We also build a holomorphic version of the product for
maps of the Riemann sphere into homogeneous spaces.  In the continuous
case we define a related module structure on the homology of maps from
a mapping cone into a manifold, and then describe a spectral sequence
that can compute it.  As a consequence we deduce a periodicity and
dichotomy theorem when the source is a compact Riemann surface and the
target is a complex projective space.
\end{abstract}

\begin{asciiabstract}
We apply a version of the Chas--Sullivan--Cohen--Jones product on the
higher loop homology of a manifold in order to compute the homology of
the spaces of continuous and holomorphic maps of the Riemann sphere
into a complex projective space.  This product makes sense on the
homology of maps from a co-H space to a manifold, and comes from a
ring spectrum.  We also build a holomorphic version of the product for
maps of the Riemann sphere into homogeneous spaces.  In the continuous
case we define a related module structure on the homology of maps from
a mapping cone into a manifold, and then describe a spectral sequence
that can compute it.  As a consequence we deduce a periodicity and
dichotomy theorem when the source is a compact Riemann surface and the
target is a complex projective space.
\end{asciiabstract}

\maketitle

\section{Introduction}

Spaces of maps, holomorphic and continuous, from a Riemann surface
into complex projective spaces
and more generally flag manifolds
have been much studied in recent years.  These
are spaces that occur naturally in physics (sigma-models, topological
field theory and gauge theory), in engineering (control theory) and
also in the
theory of minimal surfaces and harmonic maps.  Our main interest in
this paper and in its sequel \cite{paolo}
is the topology and in particular the homology of these spaces.

One main interest of the present work is that it gives the first instance
where the intersection product techniques of the newly developed theory
of ``string topology" (Chas--Sullivan and Cohen--Jones \cite{CS,CJ})
are used to give complete understanding of the homology of certain
mapping spaces. In terms of specific computations, this paper
determines the homology of maps of a two sphere into complex projective space.
It then derives important consequences for the homology of spaces of maps from
a more general closed Riemann surface. It further adapts the ``string" techniques,
mainly as developed by Cohen and Jones,
to the holomorphic setting and computes the homology of the subspace
of holomorphic maps. Without use of the intersection products introduced
originally by Chas and Sullivan, these computations would be hard
to access.

In \cite{CS}, Chas and Sullivan introduced new structures on the
homology of parameterized and unparameterized loop spaces of closed
manifolds.  One main observation there was the existence of a product
on the (regraded) homology of the free loop space of a closed manifold
which was compatible via evaluation at a basepoint with the
intersection product on the homology of the manifold in question. In
\cite{CJ}, J\,D\,S Jones and R Cohen gave a homotopy theoretic and
more rigorous definition of a similar product $\bullet$ and asserted
it was the same as Chas and Sullivan's. Shortly after and jointly with
J Yan \cite{CJY} they showed that $\bullet$ did fit naturally in the
Serre spectral sequence of the evaluation fibration to yield a
multiplicative second quadrant spectral sequence (or \textit{loop
spectral sequence}) which allows explicit computations. Further
interesting extrapolations on this spectral sequence have been
introduced recently by Le Borgne \cite{L}.

Our work exploits these ideas in various ways, it expands on them
to pass from ring to module structures and then adapts them
to the holomorphic category. It also records the fact that the results
of \cite{CJ} and hence presumably \cite{CS}
need not restrict to maps from the circle but extend naturally to maps from coassociative
co--$H$ spaces.

The theorem next and its corollary expand on the work of Cohen and Jones.
The fruitful idea in \cite{CJ} was to combine pullback diagrams of Hilbert
manifolds and Thom--Pontrjagin collapse maps.
Write $\map{}(X,Y )$ for the space
of continuous maps from $X$ into $Y$ and denote by
$ev \co  \map{}(X,Y)\lrar Y$ the evaluation map $ev(f):=f(x_0)$ where
$x_0\in X$ is a basepoint. For
$M$ a closed compact (smooth) manifold of dimension $d$, let
$TM$ be its tangent bundle. For $X$ a (based) topological space,
one can define a Thom spectrum $\map{}(X,M)^{-TM}$ as follows.
Let $\nu$ be the normal bundle of an embedding $M \to \bbr^{\ell}$, and
$ev^*(\nu)$ its pullback over $\map{}(X,M)$ via evaluation.
One then defines
$\map{}(X,M)^{-TM}$ to be the $\ell$--fold desuspension of
the suspension spectrum of the Thom space of $ev^*(\nu)$.
When $X=*$, then
$M^{-TM}$ is Spanier--Whitehead dual to
$M_+$ by a well known result of Atiyah.

All spaces are assumed to be path-connected and compactly generated.

\bth\label{firsttheo}
Suppose $S$ is a coassociative (cocommutative) co--$H$ space
coacting on a space $X$. Then

{\rm(1)}\qua $\map{}(S,M)^{-TM}$ is a (homotopy commutative) ring spectrum;

{\rm(2)}\qua $\map{}(X,M)^{-TM}$ is a module over the ring spectrum
$\map{}(S,M)^{-TM}$.

\end{theorem}

This theorem as it turns out is a consequence of a more general
statement. Let $X_1, X_2$ be based topological spaces
and suppose $M$ is a closed manifold. Then there is
a pairing of spectra
\begin{equation}\label{keypairing}
\map{}(X_1,M)^{-TM}\wedge \map{}(X_2,M)^{-TM}\lrar \map{}(X_1\vee X_2,M)^{-TM}
\end{equation}
with relevant associativity properties and commuting via the
evaluation maps
with the ring multiplication $M^{-TM}\wedge M^{-TM}\lrar M^{-TM}$.
Such a pairing exists because it is possible to construct
a tubular neighborhood of $\map{}(X_1 \vee
X_2, M)$ in $\map{}(X_1,M) \times \map{}(X_2,M)$, strictly compatible with evaluation
fibrations and living over the normal neighborhood of the diagonal
$M\hookrightarrow M\times M$ (\fullref{continuousngbd}).
Part (1) of \fullref{firsttheo} for suspensions
has also been obtained by Klein in \cite{klein} using
a different approach.

Most important for us is \fullref{firsttheo} (2) in the case when
$S$ is the $n$--sphere ($n\geq 1$), $X$ an $n$--manifold and the coaction
map $X\lrar X\vee S$ is obtained by pinching the boundary of a
small embedded disk in $X$.
The following corollary is key to our homological calculations.

\bco \label{module}
Let $N,M$ be two closed compact oriented manifolds
of dimension $n$ and $d$ respectively. Then
$\bbh_*(L^nM)$ is a ring and
$\bbh_*(\map{}(N,M))$ is
a module over $\bbh_*(L^nM)$,
where $\bbh_* = H_{*+d}$ and $L^nM = \map{}(S^n,M)$. This module structure
is compatible via the evaluation maps with the intersection product on
$\bbh_*(M)$.
\eco

Note that orientability is needed here but not in earlier statements.
The ring multiplication above (which we dub simply
\textit{intersection product}) reduces to the Cohen--Jones product for
$n=1$.  Such product has been investigated by Chataur \cite{Ch} and Hu
\cite{Hu} as well. Note also that it is possible to replace singular
homology in \fullref{module} by other multiplicative generalized
homology theories in line with work in Cohen and Godin \cite{CG} (see
\fullref{string}).  It is equally possible to expand on
\fullref{firsttheo} (1) by considering bundles over $M$ with fibrewise
multiplication hence obtaining further interesting pairings between
spectra (Gruher and Salvatore \cite{GS}).

More interestingly perhaps from the physicist's point of view is to
find out whether the module and ring structures described above can be
carried out in the holomorphic category.  Let $M$ be a complex
algebraic variety and write $\hol{}(C,M)\subset\map{}(C,M)$ for the
subspace of all holomorphic maps from a compact Riemann surface $C$
into $M$.  When $C=\bbp^1$ is the Riemann sphere, we reduce this
notation to $\hol{}(M)$. This space for $M$ compact splits into a
number of connected components depending on the geometry of $M$ and
each component is of the homotopy type of a finite dimensional
complex.  It might happen of course that $\hol{}(M)$ has no
non-constant maps.  However when $M=G/P$ is a generalized flag
manifold, $G$ being a complex Lie group and $P$ a parabolic subgroup,
$\hol{}(G/P)$ has a rich topology (Boyer, Mann, Hurtubise, Milgram and
Segal \cite{BHMM,S}).  We manage to show in \fullref{rattw} that an
intersection product similar to the Cohen--Jones product exists on the
(regraded) homology of the space $\hol{}(G/P)$.

\bth \label{flag} The Thom spectrum $\hol{}(G/P)^{-T(G/P)}$
is a homotopy commutative ring spectrum.
In particular
$\bbh_*(\hol{}(G/P)) = H_{\dim_\bbr(G/P)+*}(\hol{}(G/P))$
is a  commutative graded ring.
\end{theorem}

In proving this theorem we adapt many of the constructions of
\fullref{string} to the holomorphic category; in particular
the existence and compatibility of tubular neighborhoods.
The main geometric ingredient we use
is a construction of Gravesen and Segal \cite{G} which represents
such holomorphic maps in terms of certain configuration spaces of
\textit{principal parts} (cf \fullref{rattw}).
In the case when $G/P=\bbp^n$ we further make the
observation that the intersection products that we construct
in \fullref{module} and \fullref{flag} are compatible
via subspace inclusion (\fullref{rat}).

\bpr \label{uno}
The inclusion $\hol{}(\bbp^n)\hookrightarrow L^2 \bbp^n$ induces a
homomorphism of ring spectra
$\hol{}(\bbp^n)^{-T\bbp^n} \to (L^2 \bbp^n)^{-T\bbp^n}$.
In particular $\bbh_*( \hol{}(\bbp^n)) \to \bbh_*(L^2 \bbp^n)$ is a ring
homomorphism.
\epr

We expect the ring structure
uncovered by \fullref{flag} to be highly instrumental
in extending the homology calculations
of \cite{BHMM} for based holomorphic maps to the unbased case.
We verify that this is indeed the case for $G/P=\bbp^n$.
Write for simplicity $\hol{}(n):= \hol{}(\bbp^1,\bbp^n)$ and as
customary write $\rat{}(n)$ for the subspace of \textit{based} holomorphic
maps $\bbp^1\rightarrow\bbp^n$. If $\Omega^2(\bbp^n)$ is the
corresponding subspace of based continuous maps, then we have the
diagram of inclusions
\begin{equation}\label{allspaces}
\xymatrix{
\rat{}(n)\ar@{^{(}->}[d]\ar@{^{(}->}[r] &\Omega^2\bbp^n\ar@{^{(}->}[d]\\
 \hol{}(n)\ar@{^{(}->}[r]&L^2\bbp^n
}
\end{equation}
All connected components in this case are labeled by integers
$k$ and we write $\rat{k}\subset\hol{k}\subset L^2_{k}$ the corresponding
components of degree $k$. In the holomorphic case $k$ is
non-negative and is the algebraic degree, while in the continuous case
it takes on all values. The intersection product turns out to be
additive on degree.

In studying the Serre spectral sequence relating based and unbased
mapping spaces in \eqref{allspaces}, standard properties and
comparison methods are not enough to determine all differentials.  To
bypass this shortcoming, the basic idea has been then, using the fact
that $\bbp^n$ is closed oriented, to concoct a regraded version of
this Serre spectral sequence to get a \textit{multiplicative} spectral
sequence with respect to the intersection product $\bullet$ introduced
earlier, and hence be able to decide about differentials using the
fact that they are now derivations. This basic idea goes back to
Cohen, Jones, and Yan \cite{CJY} who worked out its technical details
for the free loop space. In \fullref{spec} we show that the techniques
of \cite{CJY} can be straightforwardly adapted to higher loop spaces
eventhough the fibers are now non-connected.  More explicitly we
construct a second quadrant spectral sequence of algebras
$E^r(S^2,\bbp^n)$ which converges to $\bbh_* (L^2\bbp^n) = H_{2n+*}
(L^2\bbp^n)$. This fact as it pleasantly turns out is sufficient to
determine entirely $H_*(L^2 \bbp^n)$ with field coefficients, and by
standard comparison arguments determine $H_* ( \hol{}(n))$ as well.
This computation is in fact at the origin of this work.

Let us recall the homology of the based holomorphic subspace as
described by Cohen, Cohen, Mann and Milgram in \cite{C2M2}.  It was
there shown that $\rat{}(n)$ is a $C_2$--space (a space with an action
of the little two-discs operad); in particular it is an $H$--space and
the ``forgetful" map $\rat{}(n)\hookrightarrow \Omega^2\bbp^n$ is a
$C_2$--map up to homotopy (and hence an $H$--map) inducing a homology
monomorphism. We make the caution that this monomorphism is special to
the genus $0$ case and does not hold for maps from positive genus
curves as is discussed by the first author in \cite{K} for example.

At the prime 2, the results of \cite{C2M2} give an algebra isomorphism
$$H_*(\rat{}(n);\bbz_2) = \bbz_2[\iota, u, Q(u), Q^2(u), \dots]$$
where the algebra structure corresponds to the Pontrjagin product, and
$\iota$ and $u$ are the images of the bottom and top homology generators of
$S^{2n-1}$ via a homotopy equivalence $S^{2n-1} \simeq \rat{1}(n)$
(cf \fullref{one}). Moreover $Q^{i}(u) := Q\cdots Q(u)\in
H_{2^{i+1}n-1}(\rat{2^i}(n);\bbz_2)$ is the
$i$--fold iterated first Dyer--Lashof
operation on $u$.
 Now the homology of the continuous mapping space
$\Omega^2 \bbp^n$ is obtained from the homology of $\rat{}(n)$ by
inverting multiplicatively $\iota$, compatibly with the fact that
$\Omega^2 \bbp^n$ is the (topological) group completion of $\rat{}(n)$.

With coefficients in $\bbz_p$, with $p$ an odd prime, the statement is
analogous, for we have
$$H_*(\rat{}(n);\bbz_p) = E[u, Q(u), Q^2(u), \dots]\tensor\bbz_p[\iota,
\beta Qu,\beta Q^2u,\dots ],
$$
where $\beta$ is the mod--$p$ homology Bockstein, while
$H_*(\Omega^2\bbp^n;\bbz_p)$ is obtained by inverting $\iota$.
Rationally the situation is trivial as all positive components of
$\rat{}(n)$ and all components of $\Omega^2 \bbp^n$ are rationally
equivalent to $S^{2n-1}$ (cf \fullref{one}). Note that
multiplication by $\iota$ or $u$ switches components up by one.

To state our main computation, write $H^*(\bbp^n )=\bbz [c]/c^{n+1}$
and grade it {\sl negatively} so that $c\in H^{-2}(\bbp^n)$.  Note
that $H_*(\Omega^2 \bbp^n) \tensor H^{*}(\bbp^n)$ is an algebra with
product induced from the Pontrjagin product in the first factor and
the cup product in the second factor.

\eject
\bth \label{main}\ \
\begin{enumerate}
\item {\sc Continuous}\qua The homology $\bbh_* (L^2 \bbp^n)$ with coefficients in a
field is (additively) the homology of the differential graded algebra $H_*(\Omega^2
\bbp^n) \tensor H^{*}(\bbp^n)$ with differential $d$ such that
$d\iota=(n+1)u c^n$ and $d$ vanishes on all other generators.
\item {\sc Holomorphic}\qua The homology $\bbh_* ( \hol{}(n))$ with coefficients in
a field is (additively) the homology of the differential graded subalgebra
$H_*(\rat{}(n)) \tensor H^{*}(\bbp^n)$ with the same differential as
above.
\end{enumerate}
\end{theorem}

\bre We actually prove more. We construct in \fullref{spec} a
multiplicative spectral sequence with $E^2$ term the algebra
$H_*(\Omega^2 \bbp^n) \tensor H^{*}(\bbp^n)$ and we show that the only
differentials are those obtained via derivations from $d$ above. In
other words $\bbh_*(L^2\bbp^n )$ as an algebra has a filtration whose
associated graded is isomorphic to the homology described in
\fullref{main} (1).  To obtain the actual multiplicative structure on
$\bbh_*(L^2\bbp^n )$ one needs to solve the extension problem at the
$E_{\infty}$--level. To that end it would be useful to work at the
chain level as in F{\'e}lix, Menichi and Thomas \cite{FMT}. Analogous
statements can be made for $\hol{}(n)$ in \fullref{main} (2).  We
conjecture that there are no extension problems and that the homology
descriptions in (1) and (2) describe respectively $\bbh_* (L^2
\bbp^n)$ and $\bbh_* ( \hol{}(n))$ as algebras.  \ere

\begin{corollary} \label{list}
Let $p$ be a prime.
\begin{enumerate}
\item {\sc Collapse}\qua If $p$ divides $(n+1)$, then
$H_*(L^2\bbp^n;\bbz_p)\cong H_*(\Omega^2 \bbp^n;\bbz_p) \tensor
H_{*}(\bbp^n;\bbz_p)$ and $H_*( \hol{}(n);\bbz_p)\cong
H_*(\rat{}(n);\bbz_p) \tensor H_*(\bbp^n;\bbz_p)$.
\item {\sc Periodicity}\qua Multiplication by $\iota^p$ induces a mod--$p$
isomorphism
$$\bbh_*(L^2_k(\bbp^n);\bbz_p ) \cong
\bbh_*(L^2_{k+p}(\bbp^n);\bbz_p )$$
If $p$ does not divide $k$, then postcomposition by a
map $\bbp^n \to \bbp^n$ of algebraic degree $k$ induces an isomorphism
$\bbh_*(L^2_1 \bbp^n;\bbz_p) \cong \bbh_*(L^2_k\bbp^n;\bbz_p)$.

\item
If $n$ is even, then the spectral sequence computing
the homology mod 2 of any component of $L^2 \bbp^n$ or any positive component
of $ \hol{}(n)$ does not collapse. The same is true if $p$ is odd
and $p$ does not divide $k(n+1)$.
If $p$ is odd and divides $k$, then the spectral sequences of $L^2_k \bbp^n$
and $ \hol{k}(n)$ mod $p$ collapse.
\end{enumerate}
\end{corollary}

Part (3) of this corollary recovers in the case $n=1$ the main result
of Havlicek \cite{Hv}.

We finally turn to maps from higher genus curves into $\bbp^n$.  In
general the homology groups of $\map{}(C,\bbp^n )$ and of its
holomorphic subspace $\hol{}(C,\bbp^n )$ are hard to compute and
little is known about them \cite{K}.  An important homological
approximation theorem due to Segal \cite{S} states that the inclusion
of $k$-th components
$$i_k : \hol{k}(C,\bbp^n)\hookrightarrow\map{k}(C,\bbp^n )$$
is a homology isomorphism up to
dimension $(k-2g)(2n-1)$, where $g$ is the genus of $C$.
In \fullref{spec} we combine the calculations of \fullref{list}
to the module structure of \fullref{module}
to gain the following important insight.

\bth \label{maintw} Let $C$ be a compact Riemann surface.
\begin{enumerate}
\item {\sc Periodicity}\qua If $p$ divides $k(n+1)$, then
$$\bbh_*(\map{i}(C,\bbp^n);\bbz_p) \cong
\bbh_*(\map{i+k}(C,\bbp^n);\bbz_p).$$
In particular if $p$ divides
$n+1$ then all components have the same homology.
\item {\sc Dichotomy}\qua The
additive homology of a component of $\map{}(C,\bbp^n)$ with
coefficients in a field is isomorphic to the homology of the component
of degree either zero or one.
\end{enumerate}
\end{theorem}

In order to obtain the homology groups
component by component, it is necessary to
find a model for these mapping spaces that allows for a much deeper insight
into their geometry. Such a model built out of configuration spaces
is constructed and studied in our sequel \cite{paolo}. One finds in particular
homology torsion of order $k(n+1)$ and $n(n+1)$.

\medskip
{\bf Acknowledgment}\qua We wish to thank Daniel Tanr\'e for his constant
support. We also thank the referees for improving with their comments
earlier versions of this paper.


\section{The degree one component and rational homology}\label{one}

As will soon be apparent, the geometry of the space of holomorphic
maps $\bbp^1\lrar\bbp^n$ is essentially reflected in the degree one
component $ \hol{1}(n)$. We fix base points so that
\begin{equation}\label{tw.on}
\rat{}(n) = \{f\co \bbp^1\lrar\bbp^n\ |\ f\ \hbox{holomorphic, and}\
f([1:0])= [1:0:..:0]\}.
\end{equation}

\ble\label{spherebundle} Let $\tau (\bbp^n)$ be the sphere bundle of
unit tangent vectors of $\bbp^n$. Then there is a fiberwise
homotopy equivalence:
$$
\xymatrix{
S^{2n-1}\ar[r]^{\simeq}\ar[d]&\rat{1}(n)\ar[d]\\
\tau\bbp^n\ar[r]^{\simeq}\ar[d]^{\pi}& \hol{1}(n)\ar[d]^{ev}\\
\bbp^n\ar[r]^{=}&\bbp^n
}$$
\ele

\proof A linear map $\bbp^1 \to \bbp^n$ sending the
base point $[1:0]$ to $[1:0:\dots:0]$ has the form $[z:x] \mapsto
[z+z_0 x:z_1 x:\dots:z_n x]$, with $z_0 \in \bbc$ and $(z_1,\dots,z_n)
\neq 0$, so that $\rat{1}(n) \cong \bbc \times (\bbc^n -0)$.  Since
$\bbp^n$ is acted on transitively by $PU(n+1)$, and the stabilizer of
the base point is $1 \times U(n)$, there is a diffeomorphism
$$ \hol{k}(n) \cong  PU(n+1) \times_{U(n)} \rat{k}(n)$$
and a smooth bundle $\rat{k}(n) \to  \hol{k}(n) \to \bbp^n$, where
$U(n)$ acts on $\rat{k}(n)$ by postcomposition.  The unit sphere
$S^{2n-1} \subset \{0\} \times (\bbc^n -0 ) \subset \rat{1}(n)$ is a
$U(n)$--invariant subspace and a deformation retract.  We identify the
tangent space to $\bbp^n$ in $[1:0: \dots:0]$ with the copy $\bbc^n
\subset \bbp^n$ embedded via $(x_1,\dots,x_n) \to [1:x_1: \dots:x_n]$,
compatibly with the action of $U(n)$.  Then $PU(n+1) \times_{U(n)}
S^{2n-1}$ is the bundle $\tau(\bbp^n)$ of unit tangent vectors with
respect to the Fubini--Study metric.  The induced map of bundles
$\tau(\bbp^n) \lrar  \hol{1}(n)$ is a homotopy equivalence, by
comparing the long exact sequences of homotopy groups.  \hfill\za
\medskip

Here is an example of a calculation that does not require
Chas--Sullivan products yet.  We denote again by $\iota \in
H_0(\rat{1}(n);\bbz)$ and $u \in H_{2n-1}(\rat{1}(n);\bbz)$ the
generators.  As discussed in the introduction
we have a product
$$H_{i_1}(\rat{k_1}(n);\bbz )\tensor\cdots\tensor
H_{i_r}(\rat{k_r}(n);\bbz )\lrar
H_{i_1+\cdots +i_r}(\rat{k_1+\cdots + k_r}(n);\bbz )
$$
and the class
$u\iota^{k-1}\in H_{2n-1}(\rat{k}(n);\bbz) \cong \bbz$ is a generator
\cite{C2M2}.

\bpr\label{kn+on}\ Let $a=[\bbp^n]\in
H_{2n}(\bbp^n;\bbz)$ be the fundamental class. Then in
the integral homology Serre spectral sequence for $\rat{k}(n) \to
\hol{k}(n) \to \bbp^n$, $d_{2n}(a) = (n+1)ku\iota^{k-1}$.  \epr

\proof In the case $k=0$ we have that $d_{2n}(a)=0$,
since the fiber is a point.
The case $k=1$ follows from \fullref{spherebundle}, since in the homology Serre spectral sequence for
 $\tau\bbp^n$, $d_{2n}(a)=(n+1)u$, as $(n+1)$ is the Euler
 characteristic of $\bbp^n$.  Now observe that precomposition with the
 $k$ degree map $\bbp^1\fract{\times k}{\lrar}\bbp^1$ given by sending
 $[z_0:z_1] \mapsto [z_0^k:z_1^k]$ takes degree one maps to degree $k$
 maps and commutes with the evaluation maps as in the diagram:
$$
\xymatrix{
\rat{1}(n)\ar[r]\ar[d]^{f_k}& \hol{1}(n )\ar[d]\ar[r]^(0.6){ev}&\bbp^n\ar[d]^{=}\\
\rat{k}(n)\ar[r]\ar[r]&\hol{k}(n )\ar[r]^(0.6){ev}&\bbp^n
}
$$
We therefore obtain a map of
Serre spectral sequences between the evaluation fibrations, and
to prove our claim it suffices by comparison of spectral sequences
to show that the fiber map $f_k$ induces multiplication by $k$ on
$H_{2n-1}\cong \bbz$.
But notice that $f_k$ extends to the $k$-th power map of the H--space
$\Omega^2\bbp^n$, the inclusion $H_*(\rat{}(n)) \to
H_*(\Omega^2\bbp^n) $ is a ring monomorphism, and the diagonal induces
$\Delta_*(u)=u \tensor \iota + \iota \tensor u$. It follows that
$f_{k*}(u)=k u \iota^{k-1}$ as desired.  \hfill\za

\bco \label{rational} For $k>0$, $H_*( \hol{k}(n);\bbq )\cong
H_*(\bbp^{n-1};\bbq)\oplus s^{2n-1}{\tilde H}_*(\bbp^{n};\bbq)$, where
$s^{2n-1}$ is the formal operator that raises degrees by $2n-1$.  \eco

\proof For $k>0$ there is a rational equivalence
$H_*(\rat{k}(n);\bbq)\cong H_*(\Omega^2_k\bbp^n;\bbq)$
\cite{C2M2}.  But $\Omega^2_k\bbp^n\simeq \Omega^2S^{2n+1}$ is
rationally the sphere $S^{2n-1}$, and its fundamental class is hit by
$[\bbp^n]$ according to \fullref{kn+on}.  \hfill\za

\section{Intersection products  and module structures}\label{string}

In this section
we adapt the intersection product of Chas--Sullivan--Cohen--Jones
to obtain homology pairings and module maps in
slightly more general settings.
Our homology pairings are not induced from maps of spaces but
are a formal consequence of constructions performed at the level
of \emph{spectra} as first indicated in \cite{CG}.
A chain level analog of the homology constructions of this section
are given in \fullref{spec}.

Throughout $M$ will refer to a closed compact $d$--manifold,
and $x_0$ is a basepoint in $S^n$.
Let $M\lrar M\times M$ be the diagonal with
normal bundle isomorphic to $TM$. The Pontrjagin--Thom collapse map
takes $M\times M$ into the Thom space $M^{TM}$. If $\tau$ is any other
bundle over $M$, $M^{\tau}$ its Thom complex, then one can
\emph{twist} the Thom--Pontrjagin construction by $\tau$
and obtain a map of based spaces
\begin{equation}\label{thompairing}
M^{\tau}\wedge M^{\tau}\lrar M^{TM \oplus 2\tau}.
\end{equation}
If $\tau$ is a virtual bundle, then the
above construction still makes sense at the level of spectra \cite{CJ}.
More precisely,
if $M\hookrightarrow \bbr^N$ is an embedding with normal
bundle $\nu$, then the Thom spectrum $M ^{-TM}$ is by definition
the $N$--fold
desuspension of the suspension spectrum of the Thom space $M^\nu$. So set $\tau = -TM$ the
virtual normal bundle of $M$. The pairing in \eqref{thompairing}
becomes
\begin{equation}\label{multiplication}
M^{-TM}\wedge M^{-TM}\lrar M^{-TM}
\end{equation}
and this endows $M^{-TM}$ with the structure of a ring
spectrum \cite{C}. It is well known that
$M^{-TM}$ is Spanier--Whitehead dual to
$M_+$, and the multiplication \eqref{multiplication} is
dual to the diagonal $M \to M \times M$.

Consider the one-point union $X\vee S$ and
write $ev\co  \map{}(\diamondsuit ,M)\lrar M$ the evaluation at the basepoint for $\diamondsuit=M$ or $\diamondsuit=S$.
We have a commutative diagram
\begin{equation}\label{moddiagram}
\xymatrix{
\map{}(X\vee S,M)\ar@{^{(}->}[r]\ar[d]^{ev_{\infty}}&
\map{}(X,M)\times \map{}(S,M)
\ar[d]^{ev\times ev}\\
M\ar@{^{(}->}[r]^{\Delta}&M\times M
}
\end{equation}
where $ev_{\infty}\co  \map{}(X\vee S, M)\lrar M$ is evaluation at
the wedgepoint and the righthand square is a pullback square.
The hooked arrows are embeddings of finite codimension
and the top such map is a map between infinite dimensional manifolds.
We wish to compare the tubular neighborhoods of top and bottom embeddings.
The following lemma is used in a special case in \cite{CJ}
but not explicitly stated nor proved there.

\begin{lemma}\label{continuousngbd}
Let $X, S$ be based spaces and $M$ a closed compact
smooth manifold.  Then the continuous mapping space $\map{}(X \vee S,
M)$ has a neighbourhood in $\map{}(X,M) \times \map{}(S,M)$
homeomorphic to the pullback of the tangent bundle of $M$ along the
evaluation map at the wedgepoint.
\end{lemma}

\proof
The point is to construct a tubular neighbourhood of the top
right embedding in diagram \eqref{moddiagram} of the form
$(ev\times ev)^{-1}(U)$ where $U$ is a tubular neighborhood of the
diagonal $\Delta \co  M\lrar M\times M$ and
with the main property that if $(DM,SM)$ is the pair
(disc,sphere)--bundles associated to the tangent bundle of $M$, then
up to homeomorphism
$$((ev\times ev)^{-1}(U),
(ev\times ev)^{-1}(\partial U))\cong (ev_{\infty}^*DM,
ev_{\infty}^*SM)$$
and these in turn are identified with
$(\tilde\nu_D,\tilde\nu_S)$, the pair (disc,sphere)--bundles associated
to the pullback bundle $\tilde\nu = ev_{\infty}^*(TM)$.

Choose a metric on $M$ with injectivity radius greater than 2
and let as above $DM\subset TM$ be the unit disc bundle of the tangent bundle
of $M$. We will write $\exp$ the exponential map $DM\lrar M$ and $\exp_x$
its restriction to $D_xM:=T_xM\cap DM$, for $x \in M$.
The image $U\subset M\times M$ of $DM$ via $(x,v) \mapsto (\exp_x(v),\exp_x(-v))$
is a tubular neighbourhood
of the submanifold $\Im(\Delta) \subset M \times M$.

We now build a map $\EXP\co  DM \to \Homeo(M)$ to the group of self
homeomorphisms of $M$ such that, for $(x,v) \in D_x(M)$,
$\EXP(x,v) =\exp_x(v)$.

For $(x,v)\in D_xM$, consider the self map $\phi_v\co  T_x(M) \to T_x(M)$
that is the identity on vectors of norm $\geq 2$ and such that
$\phi_v(w)=w + (1-|w|/2)v$ for $w \in T_x(M)$ and $|w| \leq 2$.
This map is a homeomorphism taking the origin to $v$.
The composition
$\exp_x \circ \phi_v \circ \exp_x^{-1}$ extended by the identity
outside of the exponential image of the radius 2 disc at $x$
defines the self homeomorphism $\EXP(x,v)$ of $M$.

With $ev_{\infty}$
as in \eqref{moddiagram},
the pullback $ev_{\infty}^*(DM)$ consists of pairs
$((f_1,f_2), (p,v))$ where $ev(f_1) = ev(f_2) = p$ and $v\in T_pM$,
$|v|\leq 1$. Consider then the map
$$\phi \co  ev_{\infty}^*(DM)\lrar \map{}(X,M) \times
\map{}(S,M)$$ sending
$((f_1,f_2), (p,v))\longmapsto (\EXP(p,v) \circ f_1, \EXP(p,-v) \circ f_2)$.
This map is one-to-one because the injectivity radius is greater
than 2.
The map $\phi$ is also an embedding onto
$(ev \times ev)^{-1}(U)$. Indeed, for base points
$x_0 \in X$ and $s_0 \in S$,
suppose that $g\co X \to M$ and $h\co S \to M$ are such that $(g(x_0),h(s_0)) \in U$.
 Then there exist $x \in M$ and $v \in T_x(M)$ such that $g(x_0)=\exp_x(v)$ and
$h(s_0)=\exp_x(-v)$.
Let us set $f_1 = \EXP(x,v)^{-1} \circ g$   and $f_2 = \EXP(x,-v)^{-1} \circ h$  .
By construction $\EXP(x,v)(x) = \exp_x(v)=g(x_0)$, so that
 $f_1(x_0)=x$. Similarly $f_2(s_0)=x$.  Then $(g,h)$ is $\phi((f_1,f_2),(x,v))$
 and $\phi$ is surjective onto $(ev \times ev)^{-1}(U)$ as claimed.
\hfill\za

\bre Observe that the existence of the above neighborhood depended
essentially on the fact that a smooth (connected) manifold is a
``homogeneous'' space $\Homeo(M)/\Homeo^*(M)$ of all homeomorphisms
modulo based homeomorphisms.  An analog of this construction in the
holomorphic category appears in \fullref{transverse}.  \ere

Going back to \eqref{moddiagram}, the Thom--Pontrjagin construction with respect
to the pullback $ev_{\infty}^*TM$, which we also write $TM$, yields then a map
$\map{}(X,M)\times \map{}(S,M)\lrar$
$\map{}(X\vee S,M)^{TM}$
which covers via the evaluations the map
$M\times M\lrar M^{TM}$. We can twist as before by $\tau = -TM$
and obtain the diagram of Thom spectra
\begin{equation}\label{diagramofspec}
\xymatrix{
\map{}(X,M)^{-TM}\wedge \map{}(S,M)^{-TM}
\ar[r]^(0.65){ev\wedge ev}\ar[d]&M^{-TM}\wedge M^{-TM}\ar[d]\\
\map{}(X\vee S,M)^{-TM}\ar[r]^(0.6){ev_{\infty}}&M^{-TM}\\
}
\end{equation}
We use this diagram of spectra to derive \fullref{firsttheo}
of the introduction.

\proof[Proof of \fullref{firsttheo}]
Let $S$ be a co--$H$ space and set $X=S$. Consider the composite
\begin{align*}
\mu\co  \map{}(S,M)^{-TM}\wedge \map{}(S,M)^{-TM}&\lrar \map{}(S\vee S,M)^{-TM}
\\
&\fract{m}{\lrar} \map{}(S,M)^{-TM}\end{align*}
where the first map is the left vertical of \eqref{diagramofspec} and the
second map is induced from the co--$H$ structure $S\lrar S\vee S$.
We need show that: (i) $\mu$ has a unit, (ii) is commutative if $S$ is
cocommutative and (iii) is associative if $S$ is coassociative.
As pointed out in \cite{CJ}, the unit is the composite
$S^0\lrar M^{-TM}\lrar\map{}(S,M)^{-TM}$ where the last map is
induced from the canonical section $M\lrar\map{}(S,M)$
sending $x\in M$ to the constant map at $x$.

Let's simplify the notation by setting $P:=\map{}(S,M)$,
$P\times_M\cdots\times_MP := \map{}(S\vee\cdots\vee S,M)$.
To show (ii) notice there is an involution $\tau$ on $P\times_MP=\map{}(S\vee S,M)$
which comes from permuting both factors in
$S\vee S$. When $S$ is cocommutative, $m \circ \tau\simeq m$ so that
the composite at the level of spectra
$$
P^{-TM}\wedge P^{-TM}\lrar
(P\times_MP)^{-TM}\fract{\tau^{-TM}}{\ra 4}(P\times_MP)^{-TM}
\lrar P^{-TM}
$$
is homotopic to $\mu$.
Associativity (iii) on the other hand follows
from the diagram:
\begin{equation}\label{bigdiagram}
\xymatrix{
P\times P\times P&
(P\times_MP)\times P\ar@{_{(}->}[l]\ar[r]
&P\times P\\
P\times (P\times_MP)\ar@{^{(}->}[u]\ar[d]&
P\times_MP\times_MP\ar[d]^{id\times m}\ar@{_{(}->}[l]
\ar[r]^(0.6){m \times id}\ar@{^{(}->}[u]&
P\times_MP\ar[d]^m\ar@{^{(}->}[u]\\
P\times P&P\times_MP\ar@{_{(}->}[l]
\ar[r]^m&P
}
\end{equation}
The lower-right square commutes
because $S$ is coassociative.
The other three subsquares are pullbacks, and
the hooked arrows are embeddings
with normal bundle isomorphic to the pullback of $TM$.
Then we can apply the Thom--Pontrjagin construction twisted
by $-TM$. This changes the direction of the hooked arrows
in diagram
\eqref{bigdiagram} and produces a new diagram \eqref{bigdiagram}$'$
at the level of Thom spectra.
The top left diagram of \eqref{bigdiagram}$'$ homotopy commutes
just because (a) it is obtained by taking the Thom--Pontrjagin construction
with respect to two tubular neighborhoods of
$P\times_MP\times_MP\hookrightarrow P\times P\times P$
obtained as pullbacks of tubular neighborhoods of the thin diagonal
in $M^3$, and (b) the Thom--Pontrjagin construction
doesn't depend up to homotopy on the choice of tubular neighborhood.
The three subsquares of \eqref{bigdiagram}$'$ coming from pullback squares
are hence homotopy commutative and thus the outer
square of that diagram is homotopy commutative which proves associativity.

The proof for the module structure is completely analogous. In that case
the analog of $\mu$
is obtained precomposing by the coaction map $X \to X\vee S$ and the
analog of \eqref{bigdiagram} is obtained by replacing
the rightmost $P$'s by $\map{}(X,M)$ throughout.
\hfill\za

\bre The arguments above show more generally that if
$X,Y$ and $Z$ are spaces then there is a map
\begin{align}\label{ass}
\map{}(X,M)^{-TM}\wedge \map{}(Y,M)^{-TM}&\wedge \map{}(Z,M)^{-TM}\\
&\lrar\map{}(X\vee Y\vee Z,M)^{-TM}\nonumber
\end{align}
which is well-defined up to homotopy (ie, the appropriate associativity
diagram commutes up to homotopy).
\ere

\subsection{Homology version} We now pass from ring and module
spectra to homology pairings. As in \cite{CG}, let $h_*$ be any
generalized homology theory such that
(i) the associated cohomology theory $h^*$ is multiplicative,
that is the underlying spectrum is a ring spectrum, and (ii) $M$ is $h_*$--oriented
to ensure the existence of the Thom isomorphism \cite[Chapter 13]{switzer}.
Let $\wedge$ be the intersection pairing on the homology
$h_*(M)=\bigoplus_{0\leq i\leq d}h_i(M)$ induced from the Thom isomorphism
with
respect to the diagonal embedding. This relates to the cohomology product
in $h^*$ via generalized Poincar\'e duality $pd$ according to the diagram:
$$
\xymatrix{
h_i(M)\tensor h_j(M)\ar[r]^{\wedge}\ar[d]^{pd}&h_{i+j-d}(M)\ar[d]^{pd}\\
h^{d-i}(M)\tensor h^{d-j}(M)\ar[r]^(0.55){\ \cup}&h^{2d-i-j}(M)
}$$
The following is a direct consequence of the commutativity of diagram
\eqref{moddiagram}.

\bco\label{cscoaction} Let $M$ be an $h_*$--oriented manifold of dimension $d$,
$S$ a co--$H$ space coacting on $X$.
There is a commutative diagram:
\begin{equation}\label{poincare}
\xymatrix{
h_i(\map{}(X,M))\tensor
h_j(\map{}(S,M))
\ar[d]^{ev_*}\ \ar[r] &
h_{i+j-d}(\map{}(X \vee S,M))\ar[d]^{ev_*}\\
h_i(M)\tensor
h_j(M)\ar[r]^{\wedge}&h_{i+j-d}(M)
}
\end{equation}
This gives a pairing $\bullet$ making
$h_{d+*}(\map{}(S,M))$ into a graded ring which is
commutative if $S$ is cocommutative, and
$h_{d+*}(\map{}(N,M))$ into a graded module
over $h_{d+*}(\map{}(S,M))$.
\eco

\proof
Let $m\co X \to X \vee S$ be the coaction map.
To obtain the pairing $\bullet$
apply $h_*$ to diagram \eqref{moddiagram} and use the Thom--Pontrjagin
construction to get an
Umkehr map going the ``wrong way" with the appropriate shift by
codimension of the hooked arrow (the embedding). This is depicted as follows
\begin{equation}\label{oldth.th}
\xymatrix{
h_i(\map{}(X,M)) \ot h_j(\map{}(S,M))\ar[r]\ar[dd]^{\bullet}&
h_{i+j}(\map{}(X,M)\wedge \map{}(S,M))\ar[d]\\
&h_{i+j}(\map{}(X\vee S,M)^{ev_{\infty}^*(TM)})\ar[d]^{m^*}\\
h_{i+j-d}(\map{}(X,M))
&h_{i+j}((\map{}(X,M)^{ev^*(TM)})\ar[l]_{\cong}
}
\end{equation}
where the top map is the slant pairing \cite[13.70]{switzer}
and the bottom map is the Thom isomorphism which is valid since $M$
oriented \cite[Theorem 14.6]{switzer}.
The composite in \eqref{oldth.th} is the desired pairing $\bullet$
compatible via the evaluation morphisms and according to \fullref{continuousngbd} to the intersection pairing on $h_*(M)$ as in \eqref{poincare}.

When $S$ is a co-associative co--$H$ space, the associativity of
$\bullet$ is a direct consequence of diagram \eqref{bigdiagram}.
Indeed consider the outer square in homology
(again homology is contravariant on all hooked arrows).
Starting with $x\tensor y\tensor z\in h_*(P\times P\times P)$
and going clockwise we obtain $x\bullet (y\bullet z)
\in h_{*-2d}(P)$. Going counterclockwise we obtain
$(x\bullet y)\bullet z$. A similar diagram gives the module structure.
\hfill\za

\proof[Proof of \fullref{module}]
The first part follows by setting $S=S^n$.
In this case $\bbh_*(L^n M):=H_{*+d}(L^nM)$ is a graded
\textit{commutative} ring.
The commutativity follows for
$n>1$ from the fact that $S^n$ is cocommutative and when $n=1$ from
the action of the circle on the loop space by rotations \cite{CS}.
We observe that this ring
has been constructed for $n=1$ in  \cite{CJ}
and for $n>1$  in  \cite{Hu} and in
\cite{Ch} using geometric bordism theory.

For the second part, we observe that
a closed based manifold $N$ of dimension $n$ is the cofiber of the top cell attaching map.
The associated coaction map $N\lrar N\vee S^n$ is obtained geometrically
by pinching the boundary of a small closed disc neighborhood of
the basepoint in $N$.
\hfill\za
\medskip

We will refer to the module structure map
\begin{equation}\label{coactioncase}
\bullet \co  \bbh_*(L^nM)\tensor\bbh_*(\map{}(N,M))\lrar\bbh_*(\map{}(N,M))
\end{equation}
as ``intersection module" map from now on.


\section{Holomorphic products I}\label{rat}

Let us specialize the diagram \eqref{moddiagram} to the case
$N=S^2$
\begin{equation}\label{oldth.tw}
\xymatrix{
L^2M\ar[d]^{ev}& L^2M\times_ML^2M\
\ar[l]_(0.6){*}\ar@{^{(}->}[r]\ar[d]^{ev_{\infty}}&
L^2M\times L^2M\ar[d]^{ev\times ev}\\
M&M\ar[l]_{=}\ar@{^{(}->}[r]^{\Delta}&M\times M
}
\end{equation}
In this section and the next we propose to
exhibit a commutative diagram as in \eqref{oldth.tw}
with $L^2M$ replaced by $ \hol{}(M)$
where $\hol{}(M)$ indicates the space of holomorphic maps from
$\bbp^1$ to a complex projective homogeneous $n$--manifold $M$.
As a direct consequence we obtain a ``holomorphic'' intersection product
$$H_i( \hol{}(M))\ot H_j( \hol{}(M))\lrar H_{i+j-2n}(\hol{}(M))$$
We give first the construction
for $M=\bbp^n$ since it is both simpler and more explicit than
the more general case.
The point here is to
exploit the fact that the H--space structure of the based holomorphic
maps $\rat{}(n)$ can
be defined by adding vectors of rational functions. In general this
addition is neither well-defined nor continuous on the entire space,
but its definition makes sense on a deformation retract of $\rat{}(n)
\times \rat{}(n)$.

\proof[Proof of \fullref{flag} for $G/P=\bbp^n$]
Write $\bbp^1=\bbc\cup\{\infty\}$ and let $\rat{}(n)$ be
as in \eqref{tw.on}.  If $f\in\rat{k}(n)$, then the restriction of $f$ to $\bbc$
is of the form $f(z) = [g(z):p_1(z):\cdots :p_n(z)]$, the components
being polynomials with $\deg p_i<\deg g=k$, and such that
$g,p_1,\ldots, p_n$ have no common zero. We therefore identify $f$
with its corresponding vector of rational functions
$(p_1/g,...,p_n/g)$. We call the zeroes of $g$ {\sl the poles} of $f$.

Now $PU(n+1)$ acts on $\bbp^n$ by isometries and the stabilizer of the
point $[1:0:\cdots :0]$ is a copy of $U(n)\subset PU(n+1)$. It follows
that $U(n)$ acts on $\rat{}(n)$ by postcomposition and the action is
given by matrix multiplication of $A \in U(n)$ and a vector
$f=(p_1/g,\dots,p_n/g)\in \rat{}(n)$. We remark that $Af$ and $f$ have
the same poles.

Let $\rat{\leq i}(n)$ be the $U(n)$--invariant subspace of $\rat{}(n)$
consisting of maps $f$ with poles of norm less than $i$.  The
inclusion $\rat{\leq i}(n)\subset \rat{\leq 2i}(n)$ has a
$U(n)$--equivariant homotopy inverse sending $f(z)$ to $f(2z)$, since
$U(n)$ acts by postcomposition, and the action commutes with
precomposition. This induces in turn a homotopy equivalence
$$PU(n+1)\times_{U(n)}\rat{\leq i}(n)\fract{\simeq}{\lrar}PU(n+1)
\times_{U(n)}\rat{\leq 2i}(n)$$
and in the limit we obtain a homotopy equivalence
$$PU(n+1)\times_{U(n)}\rat{\leq 1}(n) \simeq
PU(n+1)\times_{U(n)}\rat{}(n) = \hol{}(n).$$

The same argument as above when applied to translations $f(z)\mapsto
f(z-c)$ shows that $PU(n+1)\times_{U(n)}\hbox{Rat}^u(n)$ and $PU(n+1)
\times_{U(n)} \hbox{Rat}^l(n)$ are homotopy equivalent to $
\hol{}(n)$, where $\rat{}^u(n)$ and $\rat{}^l(n)$ consist of functions
having poles in some unit disc in the upper and lower half-planes
respectively. Here again $\rat{}^u(n)$ and $\rat{}^l(n)$ are
$U(n)$--invariant subspaces since the action by $U(n)$ on a given
rational function fixes the poles.
  Similarly, we have a homotopy
equivalence $PU(n+1) \times_{U(n)}(\hbox{Rat}^u(n)\times
\hbox{Rat}^l(n)) \simeq PU(n+1) \times_{U(n)}(\hbox{Rat}(n)\times
\hbox{Rat}(n)) \cong \hol{}(n)\times_{\bbp^n} \hol{}(n)$. Here $U(n)$
is acting diagonally, and
$ \hol{}(n) \times_{\bbp^n} \hol{}(n)$ is the space of pairs of
holomorphic maps which agree at the base point.

There is then a map
\begin{equation}\label{hspaceequivariance}
\hol{i}(n) \times_{\bbp^n} \hol{j}(n)\simeq
PU(n+1)\times_{U(n)}
(\hbox{Rat}_{i}^u(n)\times
\hbox{Rat}_{j}^l(n))
\fract{+}{\lrar}  \hol{i+j}(n).
\end{equation}
Here $+$ is induced
by the sum of vectors of rational functions
$\rat{i}^u(n)\times\rat{j}^l(n)\lrar$ $\rat{i+j}(n)$.
Such operation is now well defined, since we sum functions
with distinct poles, and it is $U(n)$--equivariant, since
$A(f+g)=Af+Ag$.

The twisted Thom--Pontrjagin
construction with respect to the pullback of the tangent bundle
$T\bbp^n$, seen as a normal bundle of the inclusion
$ \hol{i}(n)\times_{\bbp^n} \hol{j}(n)\hookrightarrow
 \hol{i}(n)\times \hol{j}(n)$ (see \fullref{transverse}), followed by composition pairing
 yields the ring spectrum product
 $$\hol{i}(n)^{-T\bbp^n} \wedge \hol{j}(n)^{-T\bbp^n} \to
 (\hol{i}(n)\times_{\bbp^n} \hol{j}(n))^{-T\bbp^n} \to
\hol{i+j}(n)^{-T\bbp^n}$$
The product is homotopy commutative because
$\rat{}(n)$ is homotopy commutative by a $U(n)$--equivariant homotopy.
\hfill\za

\proof[Proof of \fullref{uno}] Consider the right-half of diagram \eqref{oldth.tw} (the
pullback diagram).  There is a holomorphic version of it mapping to
\eqref{oldth.tw}
via inclusions. The Thom spaces of the corresponding
normal bundles map to each other and hence
we get a homotopy commutative diagram of spectra
$$
\xymatrix{
\hol{}(n)^{-T\bbp^n}\wedge\hol{}(n)^{-T\bbp^n}\ar[r]\ar[d]&
\left(\hol{}(n)\times_{\bbp^n}\hol{}(n)\right)^{-T\bbp^n}\ar[d]\\
(L^2\bbp^n)^{-T\bbp^n}\wedge (L^2\bbp^n)^{-T\bbp^n}\ar[r]&
\left(L^2\bbp^n\times_{\bbp^n} L^2\bbp^n\right)^{-T\bbp^n}
}$$
Composing respectively by
\ref{hspaceequivariance} and the
 loop sum to the right we end up respectively in
$\hol{}(n)^{-T\bbp^n}$ and $(L^2\bbp^n)^{-T\bbp^n}$.
We need show these compositions are compatible and to that end it suffices to show
that the following diagram commutes up to $U(n)$--equivariant homotopy
$$
\xymatrix{
\rat{i}^u(n)\times\rat{j}^l(n)\ar[r]^(0.6){+}\ar[d]&\rat{i+j}(n)\ar[d]\\
\Omega^2_i\bbp^n\times\Omega^2_j\bbp^n\ar[r]^(0.6){*}&
\Omega^2_{i+j}\bbp^n
}$$
The homotopy is given by
$H_t(f,g)= fa_t + gb_t$, where $a_t$ (respectively $b_t$)
is a homotopy
between the identity of $S^2=\bbc\cup\{\infty\}$ and the composite
$S^2 \to S^2 \vee S^2 \to S^2$
of the  pinch map and the projection onto the first
(respectively second) summand.
Here $a_0=b_0=id$, and $a_1 \co  \bbc\cup\{\infty\}\lrar\bbc\cup\{\infty\} $
(resp. $b_1$) sends the lower half-plane
(resp. upper half-plane) to $\infty$.
Note that $fa_t+gb_t$ is not holomorphic
anymore for $t>0$, however the deformations $a$ and $b$
can be chosen so that the ``poles"
 of $fa_t$ (ie, the zeros
of $p_{0,t}$ if we write $fa_t = [p_{0,t}:p_{1,t}:\ldots : p_{n,t}]$) and
those of $gb_t$ remain fixed as $t$ varies. This ensures that
$fa_t+gb_t$ is still well defined as a map into $\bbp^n$.
\endproof

\bre It is more common in the literature (eg \cite{C2M2}) to define
instead
$$\rat{}(n) = \{f\co \bbp^1\lrar\bbp^n\ |\ f\ \hbox{holomorphic, and}\
f([1:0] )= [1:1:\ldots :1]\}.$$
In this case a function $f\in\rat{}(n)$ is identified with a vector
$(f_1/g,\ldots, f_n/g)$, such that $\deg f_i=\deg g=k$ and
$f_1,\ldots, f_n,g$ are monic polynomials without common zeroes. In
this case the componentwise product of rational functions on a
deformation retract of $\rat{}(n) \times \rat{}(n)$ defines an H--space
structure, similarly as before.  It turns out that the product and the
sum give two H--space structures on $\rat{}(n)$ which are equivalent.
More generally there are $C_2$--structures associated to these
operations that are equivalent as well.  However the product is less
suitable here since it comes from the description of $\bbp^n$ as a
toric variety rather than a homogeneous space.  \ere


\section{Holomorphic products II}\label{rattw}

We now extend the product on $\hol{}(\bbp^n)$ to $\hol{}(M)\subset L^2M$
where $M$ is a
homogeneous projective variety. This is based on a beautiful extension
of the root-pole description of holomorphic maps for the case of
projective space to this more general situation as first formulated by
Gravesen and Segal, and then greatly exploited in \cite{BHMM}.

We follow the notation of \cite{BHMM} and consider a complex connected
semisimple Lie group $G$ with Lie algebra $\mathfrak{g}$. Choose a
Cartan subalgebra $\mathfrak{h}$ and the associated positive and
negative root spaces $\mathfrak{u}^+$ and $\mathfrak{u}^-$. A
parabolic subalgebra $\mathfrak{p}$ is obtained by adjoining to
$\mathfrak{h}$ all positive and some negative root subspaces.  The
corresponding subgroup $P$ is also called parabolic, and $M=G/P$ is a
generalized flag manifold, the most general form of a complex
projective homogeneous manifold.  Let us write
$\bbh_*(\hol{}(G/P))=H_{*+2d}(\hol{}(G/P))$, if $G/P$ has complex dimension
$d$.

\proof[Proof of \fullref{flag}] The complement $\mathfrak{n}^-$ of
$\mathfrak{p}$ in $\mathfrak{g}$ corresponds to a contractible group
$N \subset G$, such that the projection $N \subset G \to G/P=M$
identifies $N$ to an open set of $M$ (called the \emph{affine} part of
$M$), whose complement is the \emph{divisor at infinity}
$M_\infty$. Alternatively and under this identification $N$ is the
open cell in the Schubert decomposition for $G/P$ \cite[Section 2]{BHMM}.
If $M=\bbp^n$ then $N\cong\bbc^n$ and $M_{\infty}$ is the hyperplane
at infinity $\bbp^{n-1}$.

Let us choose as base point $[P] \in G/P$, that corresponds to the
unit $e \in N \subset G/P$.  A based holomorphic map $f\co \bbp^1 \to M$,
$f(\infty )=e$, has by definition poles at the elements of
$f^{-1}(M_\infty)$ which form a finite set.
The poles of $f$ together with
their associated ``principal parts" which we next define, determine
entirely the space $\rat{}(G/P)$ of based holomorphic maps, very much
as they do classically for meromorphic maps into $\bbp^1$.

Consider the sheaf of meromorphic functions ${\mathcal{M}}(G/P)$ on
$\bbp^1$ defined so that\break ${\mathcal{M}}(G/P)({\mathcal U})$,
$\mathcal U$ open in $\bbp^1$, consists of holomorphic functions
${\mathcal U} \to G/P$ with image not entirely contained in the
divisor at infinity $M_\infty$. Let $O(N)$ be the sheaf defined on
open sets $\mathcal U$ of $M$ by $O(N) ({\mathcal
U}):=\hol{}({\mathcal U},N)$ the space of holomorphic maps into $N$.
The action of $N$ on $G/P$ defines an action of $\hol{}({\mathcal U},
N)$ on $\hol{}({\mathcal U}, G/P)$ for any open ${\mathcal U}$ and
hence a quotient sheaf ${\mathcal{M}}(G/P)/O(N)$ which is commonly
referred to as the sheaf of \emph{principal parts} $\mathcal{PP}$. A
global section of $\mathcal{PP}$ is called a configuration of
principal parts, and corresponds to a finite number of points $z_i$ on
$\bbp^1$ (the poles) labeled by elements of their stalk
${\mathcal{PP}}_{z_i}$ \cite{jacques}. Notice that the stalk does not
depend on its location since one can translate maps from point to
point by an automorphism of $\bbp^1$.  The natural map from
$\rat{}(G/P)$ to the configuration space of principal parts without
poles at infinity is a homeomorphism by the results in \cite{BHMM}.
The illustrative example of $G/P=\bbp^n$ is discussed at the end of
this section.

This description defines on $\rat{}(G/P)$ a homotopy associative and
commutative $H$--space structure. The sum is obtained by first defining
homotopies which separate the poles of some given two functions into
disjoint discs in $\bbc$ and then taking the union of the new
configurations of principal parts (ie the homotopies change the
poles but not the elements in the stalks).  A similar procedure
defines a $C_2$--space structure.

We now turn to the unbased mapping space $\hol{}(G/P)$.
We let $P$ act as before on the right of $G$ and
on the left of $G/P$ and hence of $\rat{}(G/P)$.
The map $G \times \rat{}(G/P) \to \hol{}(G/P)$ given by
$(g,f(z)) \mapsto g \cdot f(z)$ induces a homeomorphism
$G\times_P\rat{}(G/P) \cong \hol{}(G/P)$, because
$G$ acts transitively on $G/P$ by biholomorphisms
and $P$ is the stabilizer of the base point.

  As in
\eqref{hspaceequivariance} we wish to produce a map from
$(\rat{}(G/P)\times \rat{}(G/P))\times_PG$ into $\hol{}(G/P)$. There is however
a problem in that the action of $P$ does not behave well with respect
to our principal parts description; in particular $P$ does not act on
$N$ (viewed as a subset of $G/P$).  To circumvent these issues, we
proceed as follows.

Let $U \subset P$ be a maximal compact subgroup contained in the
reductive part of $P$. More explicitly let $\mathfrak{P}^-$ be the set
of those negative roots $\alpha$ such that $\mg_\alpha \subset \mp$.
The Levi decomposition \cite{knapp} gives $\mp = \mr \oplus \mn$,
where
$$\mathfrak{r}= \mathfrak{h} \oplus \bigoplus_{\alpha \in \mP}
( \mathfrak{g}_\alpha \oplus
\mathfrak{g}_{-\alpha}),$$
$$\mn = \bigoplus_{\alpha \notin \mP,\, \alpha<0} \mg_{\alpha}
\leqno{\hbox{and}}
$$
is the radical nilpotent ideal.  Let $R$ be the Lie group of
$\mr$. This is a deformation retract of $P$.  We have that $R$
normalizes $N$ since $[\mathfrak{r},\mathfrak{n}^-]= \mathfrak{n}^- .$
Let $U \subset R$ be a maximal compact subgroup.  Then $U$ is also
maximal in $P$.  If $V$ is a maximal compact subgroup of $G$
containing $U$ then there is an identification $V/U = G/P$.

Consider the left action of $U$ on $G/P$. The subspace $N \subset G/P$
is $U$--invariant because $U$ is contained in the normalizer of
$N$. Note that the action on $N$ can be identified to the adjoint
action.  The complement of $N$ in $G/P$, the divisor at infinity, is
clearly $U$--invariant and thus the induced $U$--action on $\rat{}(G/P)$
does not change the poles.  Similarly $U$ acts on the sheaf
${\mathcal{M}}(G/P)$ by pointwise multiplication on the left. Since
the subsheaf $O(N)$ is $U$--invariant, we get an induced action of $U$
on the quotient sheaf $\mathcal{PP}$.  The action on stalks gives in
turn an action on configurations of principal parts which corresponds
to the action on $\rat{}(G/P)$ under the natural identification.  It
follows that the H--space sum of $\rat{}(G/P)$ is $U$--equivariant.

The construction of the intersection product follows now as in the case of
projective spaces, since $\hol{}(G/P) = V \times_U \rat{}(G/P)$. We however
need the Thom--Pontrjagin construction and this makes sense if a
suitable tubular neighborhood can be constructed for the subspace of
composable holomorphic maps.  This is provided by the following lemma.

\ble \label{transverse} The subspace $\hol{}(G/P)\times_{G/P} \hol{}(G/P)
\subset \hol{}(G/P)^2$ has a tubular neighbourhood homeomorphic to the
pullback of the tangent bundle along the evaluation map.  \ele

\proof The total space of the
tangent bundle of $G/P=V/U$ is $T(V/U) = V \times_U
\mathfrak{v}/\mathfrak{u}$, where $U$ acts by the adjoint
representation, and $\mathfrak{v}, \, \mathfrak{u}$ are the Lie
algebras associated to $V,\, U$.  We fix an invariant metric on $V$, that
induces a metric on $V/U$. Then
we can identify $T(V/U)$ to $V \times_U {\mathfrak{u}}^\bot \subset V
\times_U \mathfrak{v} $.  The exponential map $\exp_V\co \mathfrak{v} \to V$
can then
be used to define a map $E\co T(V/U) \to V$ by $E[g,w]=g \, 
\exp_V(w) \, g^{-1}$, for
$g \in V$ and $w \in {\mathfrak{u}}^\bot \subset \mathfrak{v}$.
The value does not depend on the chosen representative.
Let $\Phi\co V \times V/U \to V/U$ be the natural action.
Then $\Phi(E[g,w],gU)$ agrees with $\exp_{V/U}[g,w]$, where
$\exp_{V/U}\co T_{gU}(V/U) \to V/U$ is the exponential map of $V/U$.
Since $V$ acts on $V/U$ by biholomorphisms, for each $v \in V$ and
holomorphic map $z \mapsto f(z)\in V/U$
the map $z \mapsto \Phi(v,f(z))\in V/U$ is also holomorphic.

The tubular neighbourhood consists of pairs $\Phi(E[g,w] , f_1(z)),
\Phi( E[g,-w] , f_2(z))$, where $f_1,f_2\co \bbp^1 \to G/P$ are such that
$f_1(\infty)=f_2(\infty)=gU$ and $[g,w] \in T_{gU}(G/P)$
is shorter than the
injectivity radius.  This tubular neighbourhood is the preimage via
$ev \times ev$ of a tubular neighbourhood for $\Delta\co  G/P \to
(G/P)^2$.  \hfill\za

\bre For $f\in\rat{}(G/P)$, the image curve $f(\bbp^1)$ intersects
each irreducible component of the divisor $M_{\infty}$ into a finite
number of points (with multiplicities), and the vector of integers so
defined determines the connected components of $\rat{}(G/P)$.  \ere

\subsection{Example : Complex projective space}
Think of $\bbp^n$ as $SL(n+1,\bbc )/P$ where $P$ is the subgroup of
matrices of the form $\begin{pmatrix}a&B\cr 0&C\end{pmatrix}$ with $C$
of size $n\times n$ and $a\in\bbc^*$.  Then $\bbc^n$ is identified
with the subgroup of matrices $$\begin{pmatrix}1&0\cdots 0\cr
z_1&1\cdots 0\cr \vdots&\vdots\cr z_n&0\cdots 1\end{pmatrix}$$ which
acts on $SL(n+1,\bbc )/P$ by matrix multiplication on the left.  Note
that this action leaves invariant the hyperplane at infinity which
corresponds to those matrices with the top left entry zero.  If we
choose
$U=P \cap SU(n+1)$
and $V=SU(n+1)$, then we
have the identification $\bbp^n = SL(n+1,\bbc )/P = V/U$.  Given
$f\co \bbp^1\lrar SL(n+1)/P$, then away from the poles it has matrix
representation $$\begin{pmatrix}1&0\cdots 0\cr p_1/g&1\cdots 0\cr
\vdots&\vdots\cr p_n/g&0\cdots 1\end{pmatrix}$$ with reference to
\fullref{rat}.

The action of the sheaf
of holomorphic functions $h = (h_1,\ldots, h_n)$ (away from the poles)
corresponds to pointwise addition of functions with values in $N=\bbc^n$.
It is easy then to see that a principal part of a pole at $c \in \bbc$
is a $n$--tuple of elements of form
$\sum_{i=1}^{j}a_{i}(z-c)^{-i}$, seen as
an equivalence class of germs of meromorphic
functions modulo holomorphic functions.
The pointwise addition of meromorphic functions with distinct poles
keeps the principal parts of each function. This recovers the H--space
structure of \fullref{rat}.

\proof[Proof of \fullref{flag} in the general case]
 We mimic the construction for the continuous case and
observe that: (i) a diagram as in \eqref{oldth.tw} holds if we replace
the functor $L^2$ by $ \hol{}$, (ii) The Thom--Pontrjagin construction
with respect to the tubular neighborhood constructed in \fullref{transverse} yields a map $\hol{}(G/P)\times \hol{}(G/P)\lrar
(\hol{}(G/P)\times_{G/P}\hol{}(G/P))^{T(G/P)}$ which upon twisting by $-T(G/P)$
yields a multiplication on $\hol{}(G/P)^{-T(G/P)}$.
Its associativity follows by arguing as in the proof of
\fullref{firsttheo}, where
 the lower right corner of a diagram analog to diagram \eqref{bigdiagram} homotopy  commutes
because the map defining the homotopy
associativity of $\rat{}(G/P)$ is $P$--equivariant.
Finally homotopy commutativity follows because the map
defining homotopy commutativity
of  $\rat{}(G/P)$ is $P$--equivariant.
\endproof


\section{The loop spectral sequence and computations}\label{spec}

Our next objective is to determine the intersection product and module
structure for some choices of $N$ and $M$. To that end
we follow \cite{CJY} (with minor adjustments) and introduce
a graded module structure into the Serre spectral
sequence for the evaluation fibration
\begin{equation}\label{evalM}
\bmap{}(N,M)\lrar\map{}(N,M)\lrar M
\end{equation}
so that on the $E_{\infty}$ level we obtain a graded version of
the intersection module
structure constructed in \fullref{string}.
The authors in \cite{CJY} focused on the case $N=S^1$ and $M$ simply
connected. In our case
and in general $\map{}(N,M)$ is not connected
and so our construction ends up with an extra grading.

The basic case is when $N=S^m$.
The following main calculational result is key to proving \fullref{main}.
It states in brief that by regrading suitably the Serre spectral sequence for
\eqref{evalM} we end up with a multiplicative spectral sequence.

\bpr \label{spectral} Suppose that $M$ is a closed, compact, oriented,
simply connected manifold, and assume $m\geq 1$. Then there is a
spectral sequence of rings
$$\{\bbe^{r}_{p,q}(S^m,M), -\dim M\leq p\leq 0, q\geq 0 \},$$
second quadrant in $p$ and $q$, and
converging to the ring $\bbh_*(L^mM)$ in such a way as to have
$$\bbe^{2}_{p,q} \cong H^{-p}(M, H_q(\Omega^m M)).$$
The product structure on $\bbe^2$ is given by the cup product in
cohomology with coefficients in the Pontrjagin ring $H_*(\Omega^m M)$.
\epr

\proof As in \cite{CJY},
let $C_*(L^mM)$ be the singular chain complex of $L^mM$, where
$C_p(L^mM)$ is the free abelian group on
maps $\Delta^p\lrar L^mM$; $\Delta^p$ being the standard
$p$--simplex, and consider the standard filtration of this
chain complex by
$$\{0\}\hookrightarrow\cdots\hookrightarrow F_{p-1}C_*(L^mM)\hookrightarrow
F_pC_*(L^mM)\hookrightarrow\cdots \hookrightarrow C_*(L^mM)$$
where $F_p C_*(L^mM)$ is the subchain complex generated by those
$r$--simplices

$\tilde\sigma \co 
\Delta^r\lrar L^mM$ satisfying $ev\circ \tilde\sigma
= \sigma (i_0,\ldots, i_r)$ for some $\sigma \co  \Delta^q\lrar M$ and
$(i_0,\ldots, i_r)\co  \Delta^r\lrar\Delta^q$ is a simplicial map sending the
$k$-th vertex of $\Delta^r$ to the vertex $i_k$ of $\Delta^q$ with
$0\leq i_0\leq\cdots\leq i_r\leq q$), and $q\leq p$.
\begin{equation}\label{oldfi.on}
\xymatrix{
\Delta^r\ar[r]^{\tilde\sigma}\ar[d]_{(i_0,\ldots, i_r)}&L^mM\ar[d]^{ev}\cr
\Delta^q\ar[r]^{\sigma}&M}
\end{equation}
The spectral sequence arising from the above filtration is the
Serre spectral sequence for the fibration $\Omega^mM\lrar L^mM\lrar M$.
We now wish to show that the intersection product exists at the chain
level
$$cs_* \co  C_i(L^mM)\tensor C_j(L^mM)\lrar C_{i+j-d}(L^mM)$$
and is compatible with the filtrations as follows
\begin{equation}\label{oldfi.tw}
F_p(C_*(L^mM))\tensor F_q(C_*(L^mM))\lrar
F_{p+q-d}(C_*(L^mM)).
\end{equation}

Similarly as in diagram \ref{oldth.th}
$cs_*$ is the composite of several chain maps.
\begin{align*}
\times\co\  &C_*(L^mM)\tensor C_*(L^mM)\lrar C_*(L^mM\times L^mM)\\
c_1\co\ & C_*(L^mM\times L^mM)\lrar C_*((L^mM\times_M L^mM)^{\tilde\nu})\\
c_2\co\ & C_*((L^mM\times_ML^mM)^{\tilde\nu})\lrar C_*((L^mM)^{\nu})\\
w\co\  & C_*((L^mM)^{\nu})\lrar C_*(DL,SL) \\
\cap\co\  & C_*(DL,SL) \to C_{*-d}(DL).
\end{align*}
Here $(DL,SL)$ is the pullback via evaluation $L^m M \to M$
 of the pair $(DM,SM)$ corresponding to (disc,sphere) bundle
 of the tangent bundle of $M$.
 The first three maps are filtration preserving being
induced from maps of spaces and $w$ is a filtration preserving
weak equivalence.
The map $w$ has the following explicit description.
Choose a representative cochain $T_M \in C^d(DM,SM)$
for the Thom class of the tangent bundle which vanishes on
degenerate simplices, and let $T_L = ev^*(T_M)$ be the corresponding
Thom class of the pullback over $L^mM$.
Excision shows that there is a chain equivalence
$C_*(DL/SL,*)\lrar C_*(DL,SL)$. Since by construction $(L^mM)^{\nu}
= DL/SL$, the map $w$ is the composite
$$C_*((L^mM)^{\nu})= C_*(DL/SL)\lrar C_*(DL/SL,*)\lrar
C_*(DL,SL)$$

The main assertion next is then to show that capping with
$T_L$ decreases filtration by $d$ sending
$$F_p C_*(DL)\fract{\cap T_L}{\lrar} F_{p-d}C_*(DL)$$
We use the argument given in \cite[Theorem 8]{CJY} and correct
in passing a slight mistatement there.
Let $\tilde\sigma\in F_p(C_rDL )$ be in filtration $p$ so that by definition
$ev_*(\tilde\sigma ) = \sigma (i_0,\ldots, i_r)$
(as in \eqref{oldfi.on} with $L^mM$ replaced by $DL$),
$i_0\leq\cdots\leq i_r\leq q\leq p$, and
\begin{align}\label{oldfi.fo} ev_*(\tilde\sigma\cap T_L) =
ev_*(\tilde\sigma)\cap T_M& =
\sigma (i_0,\ldots, i_r)\cap T_M\\& = \pm
T_M(\sigma (i_{r-d},\ldots, i_r)) \sigma (i_0,\ldots, i_{r-d}).\nonumber
\end{align}
Since $T_M$ vanishes on degenerate simplices,
this last expression is non-zero only if
$i_{r-d}<\ldots <i_r\leq q$; that is when
$i_{r-d} \leq q-d \leq p-d$. If we write $(\sigma \rfloor_{q-d} )$
the restriction of $\sigma$ to its front $(q-d)$--face, then
\eqref{oldfi.fo} shows that for $\tilde\sigma\in F_pC_*(DL)$
$$ev_*(\tilde\sigma\cap T_L) = a
(\sigma \rfloor_{q-d} )(i_0,\ldots, i_{r-d})$$
for some constant $a$,
and hence by definition $\tilde\sigma\cap T_L\in F_{p-d}C_*(DL)$
which establishes the claim, and \eqref{oldfi.tw}.

To complete the proof we observe that since the intersection product is
filtration decreasing (by $d$), it induces a pairing at every level of
the Serre spectral sequence for $L^mM\lrar M$
$$\mu_r \co  E^r_{p,s}\tensor E^r_{q,t}\lrar E^r_{p+q-d, s+t}$$
When $M$ is simply connected, $E^2_{p,s} = H_p(M, H_s(\Omega^m M))$
and $\mu_2$ takes the form
$$H_p(M, H_s(\Omega^m M))\tensor
H_q(M, H_t(\Omega^m M))\lrar H_{p+q-d}(M, H_{s+t}(\Omega^m M))$$
which by construction corresponds to the intersection product
with coefficients in the (commutative) Pontrjagin ring $H_*(\Omega^mM)$.
Translating the spectral sequence $E^r$ to the left by $d$ as in
$\bbe^r_{-s,t}:= E^r_{d-s,t}$, $s\geq 0$,
yields a spectral sequence which now converges to the
shifted homology groups $\bbh_*(L^mM)$. The spectral sequence converges because
additively it is the Serre spectral sequence.
Via Poincar\'e duality,
the $E^2$--term becomes
$H_{d-s}(M,H_t(\Omega^m M)) \cong H^{s}(M,H_t(\Omega^m M))$ and the
intersection product in homology gets replaced by the cup product.
The spectral sequence is one of algebras with differentials behaving
as derivations.
\hfill\za

\bre Two observations are in order:

(1)\qua As the proof makes explicit, \fullref{spectral}
is valid even for a disconnected
fiber as is the case when $\Omega^mM$ has components. In this case there
is an extra grading for the components, and this grading is additive under
loop sum $\Omega^mM\times\Omega^mM\lrar\Omega^mM$.
The intersection product $\bullet$ can then be viewed as
a``link" between the various components of the $L^m M$
and this is what we exploit most.

(2)\qua Working with holomorphic mapping spaces,
it is possible as above to introduce the intersection
product directly in the spectral
sequence computing the homology of $\hol{}(n)$ and more generally
$\hol{}(G/P)$.
\ere

The proof of \fullref{spectral} used nothing special about $S^m$ aside
from the fact that it is a co--$H$ space (namely in defining the map $c_2$).
In fact a completely
analogous statement can be made after replacing throughout
$S^m$ by a coassociative co--$H$
space $S$, $\Omega^mM$ by $\bmap{}(S,M)$ and
$L^mM$ by $\map{}(S,M)$.

This analogy can be taken up further. We can consider for a space $N$
the Serre spectral sequence for \eqref{evalM}
as a second quadrant spectral sequence $\bbe (N,M)$
converging to the homology of $\map{}(N,M)$ which by Poincar\'e
duality for a simply connected $M$ can be written as
$$\bbe^2_{p,q}(N,M) = H^{-p}(M, H_q\bmap{}(N,M))\ ,\ -d\leq p\leq 0\ ,\
q\geq 0$$
The following theorem admits a proof completely analogous to that of
\fullref{spectral}.

\bth\label{comparison}
Let $S$ be a coassociative co--$H$ space
coacting on a space $N$.
Then $\bbe^r(N,M)$ is
a differential graded module over $\bbe^r(S,M)$, with action
$$\bbe^r_{p,s}(S,M)\tensor \bbe^r_{q,t}(N,M)\lrar \bbe^r_{p+q,s+t}(N,M)$$
When $r=2$,  the action corresponds to the map obtained via cup
product on $H^*(M)$ and the
$H_*(\bmap{}(S,M))$--module structure of $H_*(\bmap{}(N, M))$
induced by the coaction map $N\lrar N\vee S$.
Moreover the module structure of the $E^\infty$--terms is compatible
with the $\bbh_*(\map{}(S,M))$--module structure
of $\bbh_*(\map{}(N,M))$.
\end{theorem}

The above theorem is particularly useful to us when
$N$ is a closed oriented $m$--manifold, $S$ is the $m$--sphere
and $N\lrar N\vee S$ the pinch map.

\subsection{Applications} We are now in a position to
prove theorems \ref{main} and \ref{maintw}.

\proof[Proof of \fullref{main}]
Apply \fullref{spectral} to the case $m=2$, $M=\bbp^n$.  Write
$\bbh_i = H_{2n+i}$, and grade negatively the cohomology ring
$H^*(\bbp^n)=\bbz_2[c]/c^{n+1}$, $c\in H^{-2}(\bbp^n )$.  Recall that
$\iota\in H_0(\Omega^2_1\bbp^n)$ and $u\in H_{2n-1}(\Omega^2_1\bbp^n)$
are the generators coming from the inclusion
$S^{2n-1}\simeq\rat{1}\bbp^n\hookrightarrow \Omega^2_1\bbp^n$.

{\small
\begin{center}
\begin{picture}(300, 180)(-100, -50)
\put(00, 120){\em $E^{2n}$ term}
\put(105, 0){\vector(-1, 0){85}}
\put(105, 0){\vector(0, 1){115}}
\put(-30, 0){$H^{-*}(\bbp^n)$}
\put(115, 120){$H_*(\Omega^2\bbp^n)$}
\multiput(45, 0)(0, 40){3}{\circle*{3}}
\multiput(105, 0)(0, 40){3}{\circle*{3}}
\multiput(95, 0)(0, 40){3}{\circle*{3}}
\multiput(105, 0)(0, 40){1}{\vector(-3, 2){60}}
\multiput(60, 30)(0, 40){1}{$\times (n+1)$}
\multiput(108, 3)(0, 40){1}{$\iota$}
\multiput(108, 43)(0, 40){1}{$u$}
\multiput(108, 83)(0, 40){1}{$u^2$}
\multiput(108, 98)(0, 40){1}{$Q(u)$}
\multiput(105, 93)(0, 40){1}{\circle*{3}}
\put(30, -25){$-2n$}
\put(40, -15){$c^n$}
\put(85, -25){$-2$}
\put(92, -15){$c$}
\put(127, 40){$2n-1$}
\put(123, 40){\vector(-1, 0){14}}
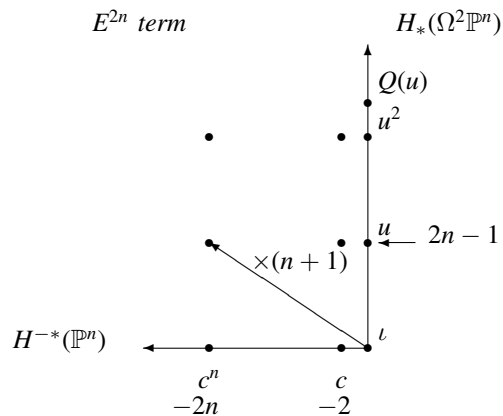
\put(-20, -40){The String Spectral Sequence Mod 2 for $L^2\bbp^n$}
\end{picture}
\end{center}}


The second quadrant spectral sequence associated to
$\Omega^2\bbp^n\lrar L^2\bbp^n\lrar\bbp^n$ (\fullref{spectral})
has as $E^2$ term the algebra $E^2_{*,*} = H^{*}(\bbp^n)\tensor
H_{*}(\Omega^2 \bbp^n)$ and converges to $\bbh_{*}(L^2\bbp^n)$.
For dimensional reasons $E^{2n+1}=E^{\infty}$.

According to \fullref{kn+on}, in the spectral sequence for
$ \hol{1}(n)$ there is a differential from the fundamental class
$[\bbp^n]$ in the base to the spherical class $u$ in the fiber. By
comparison we get the same differential for
$\Omega^2_1\bbp^n\rightarrow L^2_1\bbp^n\rightarrow\bbp^n$. Going to
the spectral sequence for $\bbh_*(L^2\bbp^n)$, the class $[\bbp^n]$
translates to $\iota\in E^2_{0,0}$, while $u$ translates to $c^n u\in
E^2_{-2n, 2n-1}$. The Euler class differential in \fullref{spherebundle} translates in turn to the differential
$d_{2n}\iota=(n+1)c^n u$. Note that the fact that $d$ is a derivation
implies that for any integer $k$
\begin{equation}\label{dtwn}
d_{2n}\iota^k = (n+1)kc^nu \iota^{k-1}
\end{equation}
This recovers the differential from \fullref{kn+on}, and also settles
the rational case (eg \fullref{rational}) since
$\Omega^2\bbp^n$ is rationally $\bbz\times S^{2n-1}$, with rational homology
$\bbq[\iota,\iota^{-1}]\tensor E[u]$.

With mod--$2$ coefficients, there are additional generators
$Q^{i}(u)$ in $H_*(\Omega^2_{2^i}(\bbp^n);\bbz_2)$.  But
$Q^{i}(u)$ is the top class in $H_*(\rat{2^i}(n);\bbz_2)$ and hence
by comparison $d_{2n}Q^{i}(u)=0$ in the spectral sequence for the
continuous mapping space.

For $p$ an odd prime, the procedure is analogous, for $Q^{i}(u)$ is
the top class in the homology mod $p$ of $\rat{p^i}(n)$, and its
Bockstein is the image of the top class of the same space with
coefficients in $\bbz_{(p)}$.

The homology of $ \hol{}(n)$ is obtained by comparison and injects
into the homology of $L^2(\bbp^n)$; the $E^{2n}$ term for $ \hol{}(n)$
being a direct summand of the $E^{2n}$ term for $L^2(\bbp^n)$. This
proves the second part of the theorem.\hfill\za
\medskip

By means of \fullref{main} one easily deduces
\fullref{list} in the introduction. We work out one example to
illustrate how the calculations go.

\bex\label{famous}
We determine all
differentials in the Serre spectral sequence $E^r$ mod 2
for the fibration
$\Omega^2 \bbp^n\rightarrow
L^2 \bbp^n\fract{ev}{\lrar}\bbp^n$.

For $n$ odd the spectral sequence collapses at the $E^2=E^{2n}$ term
because $d_{2n}$ is trivial on all multiplicative generators: $c,\iota,u$ and
$Q^i(u)$, for $i>0$.

For $n$ even and again according to \eqref{dtwn},
the differential $d_{2n}$ acts on a basis element
of the form $\iota^p u^q Q$, where $Q$ is a product of iterated Dyer--Lashof operations $Q^i(u)$, by
$d_{2n}(\iota^p u^q Q)= p\iota^{p-1}u^{q+1}Q c^n$.
Thus the surviving basis elements in the column $0$ have the form
$\iota^p u^q Q$ with $p$ even, and those in the column $-2n$ have the form
$\iota^p u^q Q c^n$ with $p$ odd.
If $$A =
\bbf_2[\iota^2,\iota^{-2},u]\tensor \bbf_2[Q^i(u),i>0] \subset H_*(\Omega^2 \bbp^n;\bbf_2)$$
then for $n$ even we have the additive isomorphism
$$H_{*+2n}(L^2\bbp^n ;\bbf_2)
= A \oplus cH_*(\Omega^2 \bbp^n;\bbf_2)\oplus \dots
\oplus c^{n-1}H_*(\Omega^2 \bbp^n;\bbf_2)  \oplus  \iota c^nA
.$$

\eex

Let now $N=C$ be a compact Riemann surface (of fixed
genus) and $M=\bbp^n$.

\proof[Proof of \fullref{maintw}] Since the components of
$\map{}(C,\bbp^n)$ are indexed by the integers, we introduce in
the spectral sequence $\bbe (C,\bbp^n)$ of \fullref{comparison}
a trigrading as in \fullref{spectral}. This spectral sequence
converges to $\bigoplus_{k\in\bbz}\bbh_*(\map{k}(C,\bbp^n))$.
The differential on $\iota^k$ in
$\bbe^{2k}_{0,0,k}(S^2,\bbp^n)$ is the class $(n+1)ku
\iota^{k-1}c^n$ which is trivial when $p$ divides $k(n+1)$. This implies
in that case that
$\iota^k$ survives to the $\bbe^{\infty}$--term and it is invertible with inverse $\iota^{-k}$.
Multiplication by $\iota^k$ switches up components by $k$ as in
$$\bbh_*(\map{i}(C,\bbp^n);\bbz_p)\fract{\iota^k}{\lrar}
\bbh_*(\map{i+k}(C,\bbp^n);\bbz_p)$$
and hence yields an isomorphism.
\hfill\za


\bibliographystyle{gtart}
\bibliography{link}

\begin{thebibliography}{}
\providecommand\bibmarginpar{\leavevmode\marginpar}
\def\urlstyle#1{{\tt #1}}

\bibitem{BHMM}
\textbf{C\,P Boyer}, \textbf{B\,M Mann}, \textbf{J\,C Hurtubise}, \textbf{R\,J
  Milgram}, \emph{The topology of the space of rational maps into generalized
  flag manifolds}, Acta Math. 173 (1994) 61--101 \xox{MR}{1294670}

\bibitem{CS}
\textbf{M Chas}, \textbf{D Sullivan}, \emph{String topology}
  \xox{arXiv}{math.GT/9911159}

\bibitem{Ch}
\textbf{D Chataur}, \emph{A bordism approach to string topology}, Int. Math.
  Res. Not.  (2005) 2829--2875 \xox{MR}{2180465}

\bibitem{C2M2}
\textbf{F\,R Cohen}, \textbf{R\,L Cohen}, \textbf{B\,M Mann}, \textbf{R\,J
  Milgram}, \emph{The topology of rational functions and divisors of surfaces},
  Acta Math. 166 (1991) 163--221 \xox{MR}{1097023}

\bibitem{C}
\textbf{R\,L Cohen},
  \href{http://projecteuclid.org/getRecord?id=euclid.hha/1139839554}
  {\emph{Multiplicative properties of {A}tiyah duality}}, Homology Homotopy
  Appl. 6 (2004) 269--281 \xox{MR}{2076004}

\bibitem{CG}
\textbf{R\,L Cohen}, \textbf{V Godin}, \emph{A polarized view of string
  topology}, from: ``Topology, geometry and quantum field theory'', London
  Math. Soc. Lecture Note Ser. 308, Cambridge Univ. Press, Cambridge (2004)
  127--154 \xox{MR}{2079373}

\bibitem{CJ}
\textbf{R\,L Cohen}, \textbf{J\,D\,S Jones},
  \href{http://dx.doi.org/10.1007/s00208-002-0362-0} {\emph{A homotopy
  theoretic realization of string topology}}, Math. Ann. 324 (2002) 773--798
  \xox{MR}{1942249}

\bibitem{CJY}
\textbf{R\,L Cohen}, \textbf{J\,D\,S Jones}, \textbf{J Yan}, \emph{The loop
  homology algebra of spheres and projective spaces}, from: ``Categorical
  decomposition techniques in algebraic topology (Isle of Skye, 2001)'', Progr.
  Math. 215, Birkh\"auser, Basel (2004)  77--92 \xox{MR}{2039760}

\bibitem{FMT}
\textbf{Y F{\'e}lix}, \textbf{L Menichi}, \textbf{J-C Thomas},
  \href{http://dx.doi.org/10.1016/j.jpaa.2004.11.004} {\emph{Gerstenhaber
  duality in {H}ochschild cohomology}}, J. Pure Appl. Algebra 199 (2005) 43--59
  \xox{MR}{2134291}

\bibitem{G}
\textbf{J Gravesen}, \emph{On the topology of spaces of holomorphic maps}, Acta
  Math. 162 (1989) 247--286 \xox{MR}{989398}

\bibitem{GS}
\textbf{K Gruher}, \textbf{P Salvatore}, \emph{String topology in a general
  fiberwise setting and applications} \xox{arXiv}{math.AT/0602210}

\bibitem{Hv}
\textbf{J\,W Havlicek}, \emph{The cohomology of holomorphic self-maps of the
  {R}iemann sphere}, Math. Z. 218 (1995) 179--190 \xox{MR}{1318152}

\bibitem{Hu}
\textbf{P Hu}, \emph{Higher string topology on general spaces}
  \xox{arXiv}{math.AT/0401081}

\bibitem{jacques}
\textbf{J\,C Hurtubise}, \emph{Configurations de particules et espaces de
  modules}, Canad. Math. Bull. 38 (1995) 66--79 \xox{MR}{1319902}

\bibitem{K}
\textbf{S Kallel}, \emph{Configuration spaces and the topology of curves in
  projective space}, from: ``Topology, geometry, and algebra: interactions and
  new directions (Stanford, CA, 1999)'', Contemp. Math. 279, Amer. Math. Soc.,
  Providence, RI (2001)  151--175 \xox{MR}{1850746}

\bibitem{paolo}
\textbf{S Kallel}, \textbf{P Salvatore}, \emph{work in progress}

\bibitem{klein}
\textbf{J\,R Klein}, \href{http://dx.doi.org/10.1090/S0002-9939-05-08148-7}
  {\emph{Fiber products, {P}oincar\'e duality and $A_\infty$--ring spectra}},
  Proc. Amer. Math. Soc. 134 (2006) 1825--1833 \xox{MR}{2207500}

\bibitem{knapp}
\textbf{A\,W Knapp}, \emph{Lie groups beyond an introduction}, Progress in
  Mathematics 140, Birkh\"auser, Boston (1996) \xox{MR}{1399083}

\bibitem{L}
\textbf{J Leborgne}, \emph{String Spectral Sequences}
  \xox{arXiv}{math.AT/0409597}

\bibitem{S}
\textbf{G Segal}, \emph{The topology of spaces of rational functions}, Acta
  Math. 143 (1979) 39--72 \xox{MR}{533892}

\bibitem{switzer}
\textbf{R\,M Switzer}, \emph{Algebraic topology---homotopy and homology},
  Classics in Mathematics, Springer, Berlin (2002) \xox{MR}{1886843}Reprint of
  the 1975 original

\end{thebibliography}

\end{document}